\documentclass[11pt,twoside]{article}

\usepackage[margin=1in]{geometry}
\usepackage{amsmath,amssymb,amsthm}
\usepackage{mathtools}
\usepackage{microtype}
\usepackage[numbers]{natbib}
\usepackage{etoolbox}
\makeatletter
\newcommand\natalpha[1]{\gdef\NAT@alpha{#1}}
\patchcmd{\NAT@bare}{\NAT@wrout{\the\c@NAT@ctr}}{\NAT@wrout{\NAT@alpha}}{}
  {\PackageError{main}{natbib patch for alphabetic labels failed}{}}
\renewcommand\NAT@spacechar{~}
\makeatother
\usepackage[shortlabels]{enumitem}
\usepackage{graphicx}
\usepackage[colorlinks=true,linkcolor=blue!60!black,citecolor=blue!60!black,urlcolor=blue!60!black]{hyperref}
\usepackage[capitalise,noabbrev,nosort]{cleveref}
\usepackage{aliascnt}
\usepackage{xcolor}
\usepackage{tcolorbox}
\usepackage{dsfont}
\usepackage[compact,small]{titlesec}
\makeatletter
\newcommand{\newaliasthm}[3]{%
  \newaliascnt{#1}{#3}%
  \let\ac@oldnewcounter\newcounter
  \renewcommand{\newcounter}[2][]{}%
  \newtheorem{#1}{#2}%
  \let\newcounter\ac@oldnewcounter
  \aliascntresetthe{#1}%
}
\makeatother
\numberwithin{equation}{section}
\numberwithin{figure}{section}

\newtheorem{theorem}{Theorem}[section]
\newaliasthm{proposition}{Proposition}{theorem}
\newaliasthm{corollary}{Corollary}{theorem}
\newaliasthm{lemma}{Lemma}{theorem}
\tcbuselibrary{skins,breakable}
\tcbset{thmbox/.style={enhanced, breakable, colback=blue!6, colframe=black,
  boxrule=0.4pt, arc=2pt, left=6pt, right=6pt, top=2pt, bottom=4pt,
  before skip=\medskipamount, after skip=\medskipamount}}
\tcolorboxenvironment{theorem}{thmbox}
\tcolorboxenvironment{proposition}{thmbox}
\tcolorboxenvironment{corollary}{thmbox}
\tcolorboxenvironment{lemma}{thmbox}
\crefname{proposition}{Proposition}{Propositions}
\crefname{corollary}{Corollary}{Corollaries}
\crefname{lemma}{Lemma}{Lemmas}
\crefformat{equation}{(#2#1#3)}
\Crefformat{equation}{(#2#1#3)}
\crefrangeformat{equation}{(#3#1#4) to~(#5#2#6)}
\crefmultiformat{equation}{(#2#1#3)}{ and~(#2#1#3)}{, (#2#1#3)}{, and~(#2#1#3)}

\newcommand{\half}{\tfrac{1}{2}}
\newcommand{\RR}{\mathbb{R}}
\newcommand{\EE}{\mathbb{E}}
\newcommand{\PP}{\mathbb{P}}
\newcommand{\ones}{\mathbf{1}}
\DeclareMathOperator{\dist}{dist}
\DeclareMathOperator{\Cov}{Cov}
\newcommand{\ip}[1]{\langle #1\rangle}
\newcommand{\indicator}[1]{\mathds{1}\{#1\}}
\newcommand{\Birkn}{\mathrm{Birk}_n}     %
\newcommand{\Cbar}{\bar{C}}              %
\newcommand{\Sn}{\mathcal{S}_n}          %
\newcommand{\rinf}{\dist_{\infty}}       %
\newcommand{\pistar}{\pi^\star}
\newcommand{\norm}[1]{\lVert #1\rVert}   %
\newcommand{\infnorm}[1]{\norm{#1}_\infty}
\newcommand{\apriori}{\emph{a priori}}
\newcommand{\Apriori}{\emph{A priori}}

\title{Smoothed analysis of stability radii for the\\ assignment problem}
\author{Michael P. Friedlander%
        \thanks{Departments of Computer Science and Mathematics,
                University of British Columbia, Vancouver, BC, Canada.
                Supported by the Natural Sciences and Engineering Research Council of Canada (NSERC),
                Discovery Grant RGPIN-2023-05250.\quad\url{https://friedlander.io},\quad \url{mpf@cs.ubc.ca}.}}
\date{October 31, 2026}

\begin{document}
\maketitle

\begin{abstract}
  Every linear program can tolerate small changes to its costs and still yield
  a solution of the original problem. The size of this tolerance measures how
  robust an optimal decision is to noise or numerical error, and how much
  regularization the problem can support. Classical sensitivity analysis 
  determines the tolerance only from an optimal solution. This paper provides
  \apriori{} lower bounds on two measures of this tolerance---the optimality
  gap and the $\infty$-norm stability radius---for the $n$-by-$n$ linear
  assignment problem over the Birkhoff polytope. The bounds derive from a
  smoothed analysis, which adds Gaussian noise with standard deviation of at
  least $\sigma$ in every direction to an arbitrary base cost. At every fixed
  confidence level, uniformly over the base cost, both quantities are at least
  of order $\sigma/(n \sqrt{\log n})$. The radius bound is asymptotically
  sharp. Under isotropic noise, the worst case occurs when every assignment has
  equal cost, and in this case, a scaled stability radius converges in
  distribution to a standard exponential variable.

  \medskip
  \noindent\textbf{Key words.} assignment problem, stability radius,
  optimality gap, smoothed analysis, sensitivity analysis, exact
  regularization, Birkhoff polytope

  \smallskip
  \noindent\textbf{MSC codes.} 90C27, 90C31, 90C05, 60D05, 68Q87
\end{abstract}

\section{\texorpdfstring{\Apriori{}}{A priori} stability}\label{sec:intro}

Every linear program can tolerate small changes to its costs and still yield a
solution of the original problem
\cite{saatyParametricObjectiveFunction1954,MangasarianMeyer1979,FriedlanderTseng2007}.
The tolerance reveals how robust an optimal decision is to noise or numerical
error, and determines how much regularization the problem can support.
Classical sensitivity analysis computes this tolerance from an optimal solution
and its reduced costs. But what can be guaranteed \emph{before} the solution is
known? A positive lower bound on the tolerance is not possible for all costs,
because the tolerance can be arbitrarily small when the two best vertices have
nearly equal cost. Even for random costs, the answer depends on the polyhedral
structure.

Here we assume that the feasible set is the Birkhoff polytope
\citep[\S II.5.11]{marcusSurveyMatrixTheory2014}, which defines the linear
assignment problem: find a minimum-cost bijection between two sets of $n$
objects.  Let $\Sn$ be the set of permutations of $\{1, \dots, n\}$. For a
cost matrix $C \in \RR^{n \times n}$, the assignment $\pi \in \Sn$ has cost
$c_\pi \coloneqq \ip{C, P_\pi}$, where $P_\pi$ is the permutation matrix with
$(P_\pi)_{i, \pi(i)} = 1$.  The assignment problem has the equivalent forms
\begin{equation}\label{eq:intro-assignment}
  \min_{\pi \in \Sn} c_\pi
  \;=\;
  \min_{X \in \Birkn} \ip{C, X},
  \qquad
  \Birkn
  \coloneqq
  \{X \in \RR_+^{n \times n} : X \ones = \ones,\: X^\top \ones = \ones\}.
\end{equation}
The two formulations have the same optimal value because the Birkhoff polytope
is the convex hull of the permutation matrices
\citep[Thm.~II.1.7]{marcusSurveyMatrixTheory2014}.

Assignment costs are often computed from noisy measurements
\citep[\S6.2.8, pp.~189--190]{burkardAssignmentProblems2009} or limited samples
\citep[\S1, Def.~1.1 and item~3]{ChenGuptaLiQiaoWang2017}. In this context, the
tolerance determines how much measurement error the optimal assignment can
absorb. We adopt the smoothed model of \citet{SpielmanTeng2004}, in which an
arbitrary deterministic base cost $\Cbar \in \RR^{n \times n}$ is perturbed to
$C = \Cbar + \sigma G$, where $G$ has independent standard Gaussian entries and
$\sigma > 0$ is the noise level.  In the applications that motivate this model,
$\Cbar$ is the ground-truth cost matrix, without measurement error. Under this
model, we establish \apriori{} probabilistic guarantees on the optimality gap
and on the $\infty$-norm stability radius, which measures how far the costs can
move before the optimal assignment changes. At every fixed confidence level,
uniformly over $\Cbar$, these quantities are at least of order
$\sigma/(n \sqrt{\log n})$. The radius bound is asymptotically sharp.

\subsection{Contributions}\label{sec:contributions}

Our main result bounds from below the probability that the optimal assignment
survives every perturbation of the costs up to a given size; the bound is
uniform in the base cost $\Cbar$ (\cref{thm:probably-exact}). The largest such
size is the $\infty$-norm \emph{stability radius} $\rinf(C)$, which measures the
distance from $C$ to the set $\Sigma$ of costs with multiple optimal
permutations (\cref{lem:nonuniqueness-distance}). For every positive threshold~$t$,
\begin{equation}\label{eq:intro-tail}
  \PP[\rinf(\Cbar + \sigma G) \ge \sigma t]
  \;\ge\;
  \exp[-O(t n \sqrt{\log n} + t^2 n)].
\end{equation}
On the event in \cref{eq:intro-tail}, the optimal permutation $\pistar$ is
unique. By the gap-radius bounds \cref{eq:cone-distance-bounds}, the
\emph{optimality gap} $\Delta(C) \coloneqq \min_{\pi \ne \pistar}(c_\pi -
c_{\pistar})$ is also bounded below: $\Delta(\Cbar + \sigma G) \ge 4 \sigma t$.

The probability bound in \cref{eq:intro-tail} is positive, uniform in $n$, when
$\sigma t$ is of order $\sigma/(n \sqrt{\log n})$. Every fixed confidence level
below one can be guaranteed by choosing a suitable threshold at this scale.
Moreover, the probability in \cref{eq:intro-tail} tends to one as $n \to
\infty$ whenever $t = o(1/(n \sqrt{\log n}))$.

The radius bound in \cref{cor:fixed-confidence} is asymptotically sharp,
uniformly over base costs, at every fixed confidence level. In particular,
additive base costs, which are constant across all assignments, minimize the
probability $\PP[\rinf(C) \ge r]$ for every $r > 0$. In this case, an
appropriately scaled stability radius converges in distribution to a standard
exponential variable (\cref{cor:extremal-base}). For a fixed noise level
$\sigma$, the same limit law holds for every sequence of base costs with
uniformly bounded entries (\cref{cor:near-additive}). We are not aware of an
existing limit law for the stability radius of a combinatorial problem under a
random or smoothed input model.

The probability bound in \cref{eq:intro-tail} is uniform in $\Cbar$, and therefore
also holds when $\Cbar$ is random and independent of $G$
(\cref{cor:random-base}). Consequently, the noise need not be isotropic, or
even Gaussian: an independent Gaussian component $\sigma G$ suffices.  This
extension includes Gaussian noise with any nonsingular covariance, with
$\sigma^2$ its smallest eigenvalue.

The stability radius and the optimality gap are two ends of one family of
distances. For $p \in [1, \infty]$, write $\dist_p(C) \coloneqq \dist_p(C,
\Sigma)$ for the \emph{distance to non-uniqueness} of $C$ in the entrywise
\mbox{$p$-norm}. Because the optimal value of \cref{eq:intro-assignment} is the
minimum of the linear functions $c_\pi$, it is differentiable in $C$ exactly
when the optimal permutation is unique. Thus $\dist_p(C)$ is also the distance
from $C$ to the set where the optimal value is not differentiable in the costs. Let $k(\pi,
\rho) \coloneqq \#\{i : \pi(i) \ne \rho(i)\}$ be the Hamming distance between
permutations $\pi$ and $\rho$. If $C$ has a unique optimal permutation
$\pistar$, then
\begin{equation}\label{eq:cone-distance}
  \dist_p(C)
  = \min_{\pi \ne \pistar} \frac{c_\pi - c_{\pistar}}{[2k(\pi, \pistar)]^{1/q}}
  \,,
  \quad\;\text{with}\quad\;
  1/p + 1/q = 1
  .
\end{equation}
This distance formula is proved in \cref{lem:nonuniqueness-distance}. At $p = 1$,
this distance is the optimality gap $\dist_1(C) = \Delta(C)$. At $p = \infty$, it
is the stability radius, as described by \citet[Thm.~3.1]{SotskovLeontevGordeev1995}.
The stability radius bounds the optimality gap as $4 \rinf(C) \le \Delta(C) \le 2n \rinf(C)$; see
\cref{eq:cone-distance-bounds}.

For $n \ge 3$ and $u > 0$, fix a permutation and the cone of costs for which it
is optimal. The costs in this cone with $\dist_p(C) \ge u$ form a translate of
the cone only when $p = \infty$ (\cref{prop:linf-translates}). The proof of
\cref{thm:probably-exact} sums the exact Gaussian measures of these translates
and therefore bounds the stability radius directly.

The event in \cref{eq:intro-tail} also makes quadratic regularization exact.
For $\mu > 0$, adding the regularizer $\half \mu \norm{X}_F^2$ to the objective
in \cref{eq:intro-assignment} defines a perturbed problem whose solution is
$P_{\pistar}$ exactly when $\mu \le 2 \rinf(C)$ (\cref{lem:margin-forms}). Thus
the exactness threshold of quadratic regularization is $2 \rinf(C)$
\cite{FriedlanderTseng2007,FKPS2025}. On the event $\{\rinf(C) \ge \sigma t\}$,
the regularized problem therefore returns the optimal assignment for every $\mu
\in (0, 2 \sigma t]$.  The upper bound $2 \sigma t$ is thus an \apriori{}
choice, available before the data is observed.

\subsection{Prior work}\label{sec:prior-work}

The isolating lemma of \citet[Lem.~1]{MVV1987} is the prototype of a gap bound
for discrete costs: if the costs are drawn independently and uniformly from
$\{1, \dots, 2n^2\}$, then the minimum-cost assignment is unique with
probability at least $\tfrac{1}{2}$, and integrality turns uniqueness into an
optimality gap of at least one.  \citet{BeierVocking2006} extend this result to
continuous weights.  For binary optimization over a feasible set in $\{0,
1\}^N$, with independent objective coefficients whose marginal densities are
bounded above by $\varphi$, their generalized isolating lemma
\citep[Lem.~5]{BeierVocking2006} establishes the bound $\PP[\Delta(C) \le s]
\le 2 \varphi N s$ for every $s \ge 0$. Vectorizing each permutation matrix places
the assignment feasible set in $\{0, 1\}^{n^2}$, so $N = n^2$. Because each
entry $C_{ij} = \Cbar_{ij} + \sigma G_{ij}$ is Gaussian with variance
$\sigma^2$, $\varphi = 1/(\sigma \sqrt{2 \pi})$. Thus, for every failure
probability $\varepsilon \in [0, 1]$, uniformly in $\Cbar$,
\begin{equation}\label{eq:beier-vocking-gap}
  \PP\Bigl[\Delta(\Cbar + \sigma G) \ge
    \varepsilon \cdot \frac{\sigma \sqrt{2 \pi}}{2 n^2}\Bigr]
  \;\ge\; 1 - \varepsilon.
\end{equation}
In contrast, at each fixed confidence level, our \cref{cor:fixed-confidence}
guarantees a gap that is larger by a factor of order $n/\sqrt{\log n}$.
Moreover, a gap guarantee of order $\sigma/(n \sqrt{\log n})$ also follows from
the anti-concentration inequality of \citet[Thm.~2.4]{BelloniFangShen2024} for
the difference of the maxima of two Gaussian vectors with arbitrary means.

Classical stability analysis characterizes the radius of an observed cost
matrix. \citet[Thm.~3.1]{SotskovLeontevGordeev1995} give the radius formula,
\cref{eq:cone-distance} at $p = \infty$, and \citet{ChakravartiWagelmans1998} characterize its
value as the root of a parameterized program. In the language of weak sharp
minima \citep{BurkeFerris1993}, $\rinf(C)$ is the largest constant $\alpha$ such
that $\ip{C, X} - c_{\pistar} \ge \alpha \norm{X - P_{\pistar}}_1$ for every $X
\in \Birkn$. At $p = \infty$, the minimum in \cref{eq:cone-distance} is attained at a vertex
of $\Birkn$ adjacent to $P_{\pistar}$. For random costs, the inequality
$\rinf(C) \ge \Delta(C)/(2n)$ turns \cref{eq:beier-vocking-gap} and the bound
derived from \citet{BelloniFangShen2024} into lower tails for the radius.  At
$\Cbar = 0$, the exact-regularization bound of \citet[\S2.1.1]{FKPS2025} is also
a lower-tail bound for the stability radius.  In our notation, their lower
bound on the mean quadratic exactness threshold \citep[Eq.~(2.2)]{FKPS2025} is
of order $\sigma N^{-3/4}$, where $N = n^2$ is the ambient dimension and the
factor $\sigma$ comes from rescaling their standard Gaussian costs. The mean is
in fact asymptotic to $\sigma/\sqrt{N \log N}$ (\cref{cor:extremal-base}). Thus
their lower bound is loose by a factor of order $N^{1/4}/\sqrt{\log N}$. The
proof of \cref{thm:probably-exact} adapts the partition into optimality cones
and the change of variables on translated cones from
\citet[Prop.~1.1, Thm.~1.3]{FKPS2025}. The theorem bounds the radius directly,
uniformly over the base cost. At each fixed failure probability, its guaranteed
radius is larger than that obtained from the Beier and Vöcking bound
\cref{eq:beier-vocking-gap} by a factor of order $n^2/\sqrt{\log n}$.

Two results on random assignment costs bear on our own. For independent uniform
$[0, 1]$ costs, which lie outside our model,
\citet[Thm.~3]{FriezeSorkin2007}
prove $\EE\,\Delta(C) \ge (1 - o(1))/n^2$.
For centered Gaussian costs, which correspond to $\Cbar = 0$ in
our model, \citet{MordantSegers2021} study the optimal value and the set of
near-optimal solutions. Our proof of \cref{thm:probably-exact} uses the same
upper bound on the expected optimal value
\citep[Prop.~2.2]{MordantSegers2021}, which determines the stability radius
scale $\sigma/(n \sqrt{\log n})$.  Their asymptotics
\citep[Thm.~2.3]{MordantSegers2021} show that this bound is sharp to first
order.  At $\Cbar = 0$, \cref{cor:extremal-base} combines their asymptotics and
the tail bounds of \cref{sec:centered-base} to establish that in the large-$n$
limit, the stability radius is exponentially distributed.

\section{Geometry of the stability radius}\label{sec:geometry}

This section collects the deterministic geometry used by our Gaussian analysis.
It describes the optimality cones, the distance to non-uniqueness at the
endpoints $p = 1$ and $p = \infty$, and equivalent forms of the radius condition
$\rinf(C) \ge u$.

\subsection{Notation and optimality cones}\label{sec:assignment-cones}

Throughout, $n \ge 2$.  We use the Frobenius inner product, and inequalities
between matrices are entrywise.  For $a, b \in \RR^n$, let $a \oplus b$ denote
the matrix with entries $(a \oplus b)_{ij} \coloneqq a_i + b_j$. We write
$\ones$ for the all-ones vector or matrix. A sequence $a_n$ is asymptotic to
$b_n$ if $a_n/b_n \to 1$ as $n \to \infty$, denoted $a_n \sim b_n$.
Probabilities $\PP$ and expectations $\EE$ are taken over all random quantities
in their argument. Unless stated otherwise, the base cost $\Cbar$ is fixed, and
all randomness comes from the noise matrix~$G$.

Recall from \cref{sec:intro} the Hamming distance $k(\rho, \pi)$, the
optimality gap $\Delta(C)$, the set $\Sigma$ of costs with multiple optimal
permutations, and the distance to non-uniqueness $\dist_p(C) = \dist_p(C,
\Sigma)$.  Any two distinct permutations differ in at least two rows.

The following identities give the norms of differences of permutation matrices,
used in \cref{lem:nonuniqueness-distance}, and the common distance from each
vertex to the barycenter $J \coloneqq \frac{1}{n} \ones \ones^\top$ of
$\Birkn$, used in \cref{lem:margin-identity}. If the Hamming distance between
$\rho$ and $\pi$ is $k$, then $\ip{P_\rho, P_\pi} = n - k$, so that
\begin{equation}\label{eq:hamming-identities}
  \norm{P_\rho - P_\pi}_F^2 = \norm{P_\rho - P_\pi}_1 = 2k,
  \quad\text{and}\quad
  \norm{P_\pi - J}_F^2 = n - 1.
\end{equation}
Recall that cost differences satisfy $c_\rho - c_\pi = \ip{C, \, P_\rho -
P_\pi}$.

For each permutation $\pi \in \Sn$, the \emph{optimality cone} of $\pi$, also
known as the inner normal cone of $\Birkn$ at $P_\pi$, is the set of costs that
$\pi$ minimizes,
\[
  \mathcal{N}_\pi \coloneqq \{C :
    \ip{C, \, P_\rho - P_\pi} \ge 0
    \:\text{for all}\: \rho \in \Sn\}.
\]
By \citet[Prop.~1.1]{FKPS2025}, specialized to $\Birkn$, the cones cover
$\RR^{n \times n}$ and overlap only on the measure-zero hyperplanes with tied
costs: $\{C : c_\rho = c_\pi\}$ with $\rho \ne \pi$.

\subsection{Distance to non-uniqueness and the stability radius}

For a cost $C$ with a unique optimal permutation $\pi^\star$, the classical
\emph{stability radius} of the optimum $\pi^\star$ in the entrywise norm
$\infnorm{E} \coloneqq \max_{i,j} \lvert E_{ij}\rvert$ is
\[
  \sup\{r \ge 0 :
    \pi^\star \text{ is optimal for } C + E
    \text{ whenever } \infnorm{E} \le r\}.
\]
The following establishes the distance formula
\cref{eq:cone-distance} and, at $p = \infty$, identifies the stability radius
with $\rinf(C)$.

\begin{lemma}[Distance to non-uniqueness]\label{lem:nonuniqueness-distance}
  Let $C$ have unique optimal permutation $\pistar$, and let $1 \le p, q \le
  \infty$ satisfy $1/p + 1/q = 1$.  In the entrywise $p$-norm,
  \begin{equation}\label{eq:dist-p-def}
    \dist_p(C)
    = \dist_p(C, \partial \mathcal{N}_{\pistar})
    = \min_{\rho \ne \pistar} \frac{c_\rho - c_{\pistar}}{(2k(\rho, \pistar))^{1/q}} \;.
  \end{equation}
  Moreover, $\dist_p(C)$ is the largest $r$ for which $\pistar$ is optimal for
  every cost $C + E$ with $\norm{E}_p \le r$, and every cost $C + E$ with
  $\norm{E}_p < \dist_p(C)$ has unique optimum $\pistar$.
\end{lemma}

The Euclidean case ($p = 2$) in \cref{lem:nonuniqueness-distance} is given by
\citet[Thm.~3.3]{SotskovLeontevGordeev1995}.

\begin{proof}[Proof of \cref{lem:nonuniqueness-distance}]
  Let $r$ be the minimum on the right.  For $\rho \ne \pistar$, set $a_\rho
  \coloneqq P_\rho - P_{\pistar}$ and $g_\rho \coloneqq c_\rho - c_{\pistar} >
  0$, so that $g_\rho = \ip{C, a_\rho}$.  The matrix $a_\rho$ has $2k(\rho,
  \pistar)$ entries equal to $\pm 1$ and all others zero, so $\norm{a_\rho}_q =
  (2k(\rho, \pistar))^{1/q}$.  If $\norm{E}_p < r$, then Hölder's inequality
  yields, for every $\rho \ne \pistar$,
  \[
    \ip{C + E, a_\rho}
    \;\ge\; g_\rho - \norm{E}_p \norm{a_\rho}_q
    \;>\; 0,
  \]
  so $C + E$ has unique optimum $\pistar$ and lies in the interior of
  $\mathcal{N}_{\pistar}$.  The open ball of radius $r$ around $C$ therefore
  misses both $\Sigma$ and $\partial \mathcal{N}_{\pistar}$, so both distances
  are at least $r$.  If $\norm{E}_p \le r$, the same bound implies $\ip{C + E,
  a_\rho} \ge 0$, so $\pistar$ is optimal for $C + E$. Conversely, let $\rho_0$
  attain the minimum, set $s \coloneqq 2k(\rho_0, \pistar)$, and set $E_0
  \coloneqq -(g_{\rho_0}/s) \, a_{\rho_0}$. Then $\norm{E_0}_p = g_{\rho_0}
  s^{1/p}/s = r$ and, by \cref{eq:hamming-identities}, $\ip{C + E_0, a_{\rho_0}}
  = g_{\rho_0} - g_{\rho_0} = 0$, while Hölder's inequality and the definition of
  $r$ give $\ip{C + E_0, a_\rho} \ge g_\rho - r \norm{a_\rho}_q \ge 0$ for every
  $\rho \ne \pistar$.  Thus $\pistar$ and $\rho_0$ are both optimal for $C +
  E_0$, which lies in $\Sigma$ and, being in $\mathcal{N}_{\pistar}$ with a
  tight inequality, on $\partial \mathcal{N}_{\pistar}$.  Both distances are
  therefore at most $r$, and hence equal to $r$.  Finally, for every $\eta > 0$,
  $\ip{C + (1 + \eta) E_0, a_{\rho_0}} = -\eta g_{\rho_0} < 0$, so $\pistar$ is
  not optimal for $C + (1 + \eta) E_0$, whose perturbation has norm $(1 + \eta)
  r$.  Hence $r$ is the largest radius for which $\pistar$ stays optimal under
  every perturbation of norm at most $r$.
\end{proof}

The two endpoints of the distance family are related by the bounds
\begin{equation}\label{eq:cone-distance-bounds}
  4 \rinf(C) \le \Delta(C) \le 2n \rinf(C).
\end{equation}
The lower bound follows because $k(\rho, \pi^\star) \ge 2$ for every non-optimal
assignment $\rho$; the upper bound follows because a competitor that attains the
minimum in \cref{eq:dist-p-def} at $p = \infty$ has excess cost at most $2n
\rinf(C)$. Direct comparison of the entrywise $1$- and $\infty$-norms yields
only $\rinf(C) \le \Delta(C) \le n^2 \rinf(C)$.

\subsection{Translated cones and quadratic exactness}

In \cref{sec:access} we compute the exact Gaussian probability of the event
$\{\rinf(C) \ge u\}$ by a change of variables, which requires decomposing the
event into disjoint translates of the optimality cones~$\mathcal{N}_\pi$.
\Cref{lem:margin-forms} provides this decomposition. These same translates
consist of the costs at which the quadratically regularized problem
\begin{equation}\label{eq:quadratic-regularization}
  \min_{X \in \Birkn} \, \ip{C, X} + \half \mu \norm{X}_F^2
\end{equation}
is exact for a range of values of the positive regularization parameter $\mu$.

\begin{lemma}[Radius regions and quadratic exactness]\label{lem:margin-forms}
  Let $u > 0$, $\pi \in \Sn$, and $C \in \RR^{n \times n}$.  The following
  conditions are equivalent:
  \begin{enumerate}[(a)]
    \item $\pi$ is optimal for $C$, and $\rinf(C) \ge u$;
    \item $C \in \mathcal{N}_\pi + 2u(J - P_\pi)$;
    \item $P_\pi$ is the solution of \cref{eq:quadratic-regularization} with
      $\mu = 2u$.
  \end{enumerate}
  In that case $\pi$ is the unique optimal permutation of $C$.
\end{lemma}

\begin{proof}[Proof of \cref{lem:margin-forms}]
  Set $w_\pi \coloneqq 2u(J - P_\pi)$.  The identities $\ip{J, P_\rho - P_\pi} =
  0$ and $\ip{P_\rho, P_\pi} = n - k(\rho, \pi)$ give $\ip{w_\pi, P_\rho - P_\pi}
  = 2u k(\rho, \pi)$. Assume (b), so that $C = A + w_\pi$ with $A \in
  \mathcal{N}_\pi$.  For every $\rho \ne \pi$,
  \[
    \ip{C, P_\rho - P_\pi} = \ip{A, P_\rho - P_\pi} + 2u k(\rho, \pi)
    \ge 2u k(\rho, \pi) > 0,
  \]
  so $\pi$ is the unique optimum and $\rinf(C) \ge u$ by \cref{eq:dist-p-def}
  at $p = \infty$, which proves (a) and the final assertion. Conversely, assume
  (a).  Since $\rinf(C) \ge u > 0$, the cost $C$ is not in $\Sigma$, so $\pi$ is
  its unique optimum.  By \cref{eq:dist-p-def} at $p = \infty$, every $\rho \ne
  \pi$ satisfies $c_\rho - c_\pi \ge 2 k(\rho, \pi) \rinf(C) \ge 2u k(\rho,
  \pi)$, so
  \[
    \ip{C - w_\pi, P_\rho - P_\pi} = c_\rho - c_\pi - 2u k(\rho, \pi) \ge 0,
  \]
  proving (b).  The gradient of the objective in
  \cref{eq:quadratic-regularization} with $\mu = 2u$ at $P_\pi$ is $C + 2u
  P_\pi$, so by convexity $P_\pi$ is the solution exactly when $C + 2u P_\pi \in
  \mathcal{N}_\pi$.  Adding a matrix $a \oplus b$ to a cost changes no cost
  difference, so $\mathcal{N}_\pi + a \oplus b = \mathcal{N}_\pi$. With $a
  \oplus b = 2u J$, the condition becomes $C \in \mathcal{N}_\pi + 2u(J -
  P_\pi)$, which is (b).
\end{proof}

The proof of \cref{thm:probably-exact} needs one further property: the set
\[
  K_{\pi, p}(u) \coloneqq \{C \in \mathcal{N}_\pi : \dist_p(C) \ge u\}
\]
of costs at which $\pi$ survives every perturbation of size up to $u$ must be a
translate of $\mathcal{N}_\pi$. \Cref{prop:linf-translates} characterizes the
entrywise $\infty$-norm as the unique entrywise $p$-norm with this translation
property.

\begin{proposition}[Translation property of the $\infty$-norm]\label{prop:linf-translates}
  Let $n \ge 2$, $\pi \in \Sn$, $u > 0$, and let $1 \le p, q \le \infty$
  satisfy $1/p + 1/q = 1$. Then $K_{\pi, \infty}(u) = \mathcal{N}_\pi + 2u(J -
  P_\pi)$. For $n \ge 3$, $K_{\pi, p}(u)$ is a translate of $\mathcal{N}_\pi$
  only if $p = \infty$. For $n = 2$, it is a translate for every $p$.
\end{proposition}

\begin{proof}[Proof of \cref{prop:linf-translates}]
  At $p = \infty$, \cref{lem:margin-forms} establishes the stated translation.
  For every $\rho \ne \pi$, set $a_\rho \coloneqq P_\rho - P_\pi$.  Its
  entrywise $q$-norm is $\norm{a_\rho}_q = (2k(\rho, \pi))^{1/q}$.  Because $u
  > 0$, every cost in either set below has unique optimum $\pi$, so by
  \cref{lem:nonuniqueness-distance},
  \begin{equation}\label{eq:distance-superlevel}
    K_{\pi, p}(u)
    = \{C : \ip{C, a_\rho} \ge u \norm{a_\rho}_q
      \:\text{for every}\: \rho \ne \pi\}.
  \end{equation}

  Suppose that $K_{\pi, p}(u) = \mathcal{N}_\pi + w$.  Because $0 \in
  \mathcal{N}_\pi$ and $\ip{A, a_\rho} \ge 0$ for every $A \in
  \mathcal{N}_\pi$, the minimum of $\ip{C, a_\rho}$ over $\mathcal{N}_\pi + w$
  is $\ip{w, a_\rho}$. Let $P_\rho$ be a vertex of $\Birkn$ adjacent to
  $P_\pi$, that is, joined to it by an edge.  The edge $[P_\pi, P_\rho]$ is a
  face of $\Birkn$, so some cost $B$ attains its minimum over $\Birkn$ exactly
  on it. Thus, $\ip{B, a_\rho} = 0$ and $\ip{B, a_\tau} > 0$ for every $\tau
  \ne \pi, \rho$. By \cref{eq:hamming-identities}, the cost $\lambda B + u
  \norm{a_\rho}_q \, a_\rho/(2k(\rho, \pi))$ attains the bound for $\rho$ in
  \cref{eq:distance-superlevel}, and for all sufficiently large $\lambda$ it
  satisfies the finitely many other bounds, in which $\lambda$ has a positive
  coefficient.  Hence
  \begin{equation}\label{eq:facet-shift}
    \ip{w, a_\rho} = u \norm{a_\rho}_q
    \quad\text{whenever } P_\rho \text{ is adjacent to } P_\pi.
  \end{equation}

  Now let $n \ge 3$, fix three rows, and let $\rho$ range over the five
  permutations other than $\pi$ ($3! - 1 = 5$) that agree with $\pi$ outside
  these rows. Each $P_\rho$ is adjacent to $P_\pi$.  Indeed, the matrices in
  $\Birkn$ supported on the union of the supports of $P_\pi$ and $P_\rho$ form
  a face: they minimize over $\Birkn$ the cost that is zero on this union and
  one elsewhere.  In each row, a permutation matrix in this face takes the
  entry of $P_\pi$ or of $P_\rho$.  So the rows where it differs from $P_\pi$
  and the rows where it differs from $P_\rho$ are disjoint subsets of the at
  most three rows where $\pi$ and $\rho$ differ.  For a permutation matrix
  other than $P_\pi$ and $P_\rho$, each subset has at least two rows, which is
  impossible.  So the face is the edge $[P_\pi, P_\rho]$.

  Three of these permutations differ from $\pi$ in two rows, and the other two
  in three rows.  In the selected rows and the columns assigned to them by
  $\pi$, the first three permutation matrices sum to the all-ones matrix, as do
  $P_\pi$ and the other two; outside this block, both sums equal $3 P_\pi$.
  Subtracting $3 P_\pi$ shows that the three directions $a_\rho$ of the first
  group sum to the two of the second.  Applying $\ip{w, \cdot}$ and
  \cref{eq:facet-shift}, then canceling $u > 0$, yields $3 \cdot 4^{1/q} = 2
  \cdot 6^{1/q}$.  Hence $(3/2)^{1/q} = 3/2$, so $q = 1$ and $p = \infty$.

  When $n = 2$, the only permutation $\rho \ne \pi$ has $\ip{a_\rho, a_\rho} =
  4$, so \cref{eq:distance-superlevel} reduces to $K_{\pi, p}(u) =
  \mathcal{N}_\pi + (u 4^{1/q}/4) \, a_\rho$ for every $p$.
\end{proof}

\section{A lower tail for the stability radius}\label{sec:access}

This section proves \cref{thm:probably-exact}, a lower-tail bound for the
stability radius that is uniform over the base cost. By \cref{lem:margin-forms},
the radius event $\{\rinf(C) \ge u\}$ is a disjoint union of translated
optimality cones. We use a change of variables to rewrite the Gaussian measure
of each translate as an expectation over the original cone $\mathcal{N}_\pi$.
Because the cones tile the space and overlap only on their boundaries, these
expectations sum to a single expectation (\cref{lem:margin-identity}). The base
cost $\Cbar$ enters this expectation only through $\pi^\star(C)$, the almost
surely unique optimal permutation. Replacing the integrand by its minimum over
all permutations yields a lower bound independent of $\Cbar$.
\Cref{cor:random-base} extends the theorem to random base costs and to
anisotropic Gaussian noise.

\begin{theorem}[Uniform lower tail for the stability radius]\label{thm:probably-exact}
  Let $n \ge 2$, and let $\Cbar \in \RR^{n \times n}$ and $\sigma > 0$ be
  arbitrary. Set $C = \Cbar + \sigma G$, where $G$ has i.i.d.\ standard
  Gaussian entries.  Then for every $t > 0$,
  \begin{equation}\label{eq:gap-law}
    \PP[\rinf(C) \ge \sigma t]
    \;\ge\;
    \exp\Bigl(
      -2 t \sqrt{2(n-1) \log(n!)}
      \;-\; 2 (n-1) t^2
    \Bigr).
  \end{equation}
  On this event the optimal permutation $\pi^\star$ is unique,
  \begin{equation}\label{eq:margin-inequality}
    c_\rho - c_{\pi^\star} \ge 2 \sigma t \, k(\rho, \pi^\star)
    \quad\text{for every}\quad \rho \ne \pi^\star,
  \end{equation}
  and $\Delta(C) \ge 4 \sigma t$.
\end{theorem}

\begin{lemma}[Gaussian radius identity]\label{lem:margin-identity}
  In the setting of \cref{thm:probably-exact}, for every $u \ge 0$,
  \begin{equation}\label{eq:smoothed-mixture}
    \PP[\rinf(C) \ge u]
    \;=\; e^{-2(n-1) u^2/\sigma^2} \,
    \EE \exp\Bigl(-\frac{2u}{\sigma} \ip{G, J - P_{\pi^\star(C)}}\Bigr).
  \end{equation}
\end{lemma}

\begin{proof}[Proof of \cref{lem:margin-identity}]
  The identity at $u = 0$ is immediate.  Fix $u > 0$. For each $\pi$, translate
  the optimality cone $\mathcal{N}_\pi$ by the vector
  \begin{equation}\label{eq:shift-vector}
    w_\pi \coloneqq 2u (J - P_\pi),
    \qquad
    \norm{w_\pi}_F^2
    = 4u^2 \, \norm{J - P_\pi}_F^2
    = 4u^2 (n-1),
  \end{equation}
  where the norm follows from \cref{eq:hamming-identities}. This shift has
  minimum Frobenius norm among all shifts producing the same translate: the
  vector $w_\pi$ is orthogonal to every matrix $a \oplus b$, and adding such a
  matrix leaves each cone unchanged. By \cref{lem:margin-forms}, the radius
  event is the disjoint union
  \[
    \{C : \rinf(C) \ge u\} = \bigcup_{\pi \in \Sn} \, (\mathcal{N}_\pi + w_\pi).
  \]
  Each translated cone consists of costs with unique optimum $\pi$.

  Write $\gamma_{\Cbar, \sigma}$ for the law of $C$ and $\varphi_{\Cbar, \sigma}$
  for its density.  Because the translated cones are disjoint, $\PP[\rinf(C)
  \ge u] = \sum_{\pi \in \Sn} \gamma_{\Cbar, \sigma}(\mathcal{N}_\pi + w_\pi)$.
  For each term, substitute $C = y + w_\pi$ and expand the square in the
  exponent of the density:
  \begin{equation}\label{eq:shifted-cone-measure}
    \begin{aligned}
      \gamma_{\Cbar, \sigma}(\mathcal{N}_\pi + w_\pi)
      &= \int_{\mathcal{N}_\pi} \varphi_{\Cbar, \sigma}(y + w_\pi) \, dy \\
      &= (2 \pi \sigma^2)^{-n^2/2} \int_{\mathcal{N}_\pi}
        \exp\biggl(
          -\frac{\norm{y - \Cbar}_F^2 + 2 \ip{y - \Cbar, \, w_\pi}
            + \norm{w_\pi}_F^2}{2 \sigma^2}
        \biggr) \, dy \\
      &= e^{-\norm{w_\pi}_F^2/(2 \sigma^2)}
        \int_{\mathcal{N}_\pi} e^{-\ip{y - \Cbar, \, w_\pi}/\sigma^2} \,
        \varphi_{\Cbar, \sigma}(y) \, dy \\
      &= e^{-2(n-1) u^2/\sigma^2} \,
        \EE\bigl[\indicator{\Cbar + \sigma G \in \mathcal{N}_\pi} \,
          e^{-\ip{G, \, w_\pi/\sigma}}\bigr],
    \end{aligned}
  \end{equation}
  where the integral is rewritten as an expectation over $G$ using $y = \Cbar +
  \sigma G$ and the norm identity in \cref{eq:shift-vector}. The cones
  $\{\mathcal{N}_\pi\}_{\pi \in \Sn}$ cover $\RR^{n \times n}$ and overlap only
  on finitely many hyperplanes (\cref{sec:assignment-cones}), which have
  probability zero under the Gaussian law of $C$. Hence summing
  \cref{eq:shifted-cone-measure} over $\pi$ yields \cref{eq:smoothed-mixture},
  because $w_\pi = 2u (J - P_\pi)$.
\end{proof}

\begin{proof}[Proof of \cref{thm:probably-exact}]
  Fix $t > 0$.  The deterministic assertions hold on the event $\{\rinf(C) \ge
  \sigma t\}$: there $\rinf(C)$ is positive, so $C \notin \Sigma$ and
  $\pi^\star$ is unique.  Moreover, \cref{eq:dist-p-def} at $p = \infty$
  implies $c_\rho - c_{\pi^\star} \ge 2 k(\rho, \pi^\star) \rinf(C)$ for every
  $\rho \ne \pi^\star$, which establishes \cref{eq:margin-inequality}. Since
  $k(\rho, \pi^\star) \ge 2$, we also have $\Delta(C) \ge 4 \sigma t$.  It
  remains to bound $\PP[\rinf(C) \ge \sigma t]$.

  Apply \cref{lem:margin-identity} with $u = \sigma t$, so that its prefactor
  is $e^{-2(n-1) t^2}$, and set $H_\pi \coloneqq \ip{G, \, P_\pi - J}$.
  Jensen's inequality and the pointwise bound
  \[
    H_{\pi^\star(C)} \ge \min_\pi H_\pi
  \]
  give
  \begin{equation}\label{eq:jensen}
    \EE \, \exp(2 t H_{\pi^\star(C)})
    \;\ge\; \exp(2 t \, \EE H_{\pi^\star(C)})
    \;\ge\; \exp(-2 t \, \EE \max_\pi (-H_\pi)).
  \end{equation}
  Each $-H_\pi$ is centered Gaussian with variance $\norm{P_\pi - J}_F^2 = n -
  1$. The expected maximum of $M$ centered Gaussian variables, each with
  variance at most $v$, is at most $\sqrt{2 v \log M}$, whether or not the
  variables are independent \citep[Thm.~2.5]{BoucheronLugosiMassart2013}.
  Applied with $M = n!$ and $v = n - 1$ to the expected maximum on the right of
  \cref{eq:jensen}, this bound gives
  \begin{equation}\label{eq:Bn}
    \EE \max_\pi (-H_\pi) \;\le\; \sqrt{2(n-1) \log(n!)} \, .
  \end{equation}
  Combining this bound with \cref{eq:jensen} and \cref{eq:smoothed-mixture}
  proves~\cref{eq:gap-law}.
\end{proof}

The quantities in~\cref{eq:Bn} reappear later, so we give each its own symbol:
\begin{equation}\label{eq:Yn-Bn}
  Y_n \coloneqq \max_{\pi \in \Sn} \, \ip{G, J - P_\pi},
  \qquad
  B_n \coloneqq \sqrt{2(n-1) \log(n!)} \, .
\end{equation}
Proposition~2.2 of \citet{MordantSegers2021} also establishes the
bound~\cref{eq:Bn}, and their Theorem~2.3 establishes $\EE Y_n/B_n \to 1$.
Thus~\cref{eq:Bn} is sharp to first order. The bound $B_n$ sets the stability
radius scale $\sigma/(2B_n)$ used by \cref{cor:fixed-confidence,cor:extremal-base},
and $Y_n$ determines the radius tail in \cref{sec:centered-base}. By Stirling's
approximation, $B_n \sim n \sqrt{2 \log n}$, so the scale $\sigma/(2B_n)$ is of
order $\sigma/(n \sqrt{\log n})$ as stated in~\cref{sec:intro}.

The change of variables in \cref{eq:shifted-cone-measure} and the Jensen bound
in \cref{eq:jensen} adapt the corresponding steps in the proof of
\citet[Thm.~1.3]{FKPS2025}, which bounds the centered Gaussian measure of each
translated cone. Here, however, the law of $C$ may be noncentered. Thus the
measure of each translated cone depends on $\Cbar$, and only their sum has a
positive lower bound uniform in $\Cbar$. The proof therefore sums the Gaussian
measures of the translated cones before applying Jensen's inequality. The two
proofs differ in how they bound $\EE H_{\pi^\star(C)}$ in \cref{eq:jensen}.
\citet{FKPS2025} use the Cauchy--Schwarz inequality, which yields, in our
notation, a linear term of order $t \, n^{3/2}$ in the exponent at $\Cbar = 0$
\citep[\S2.1.1]{FKPS2025}. The maximal inequality over the $n!$ vertices of
$\Birkn$ in \cref{eq:Bn} yields the linear term $2 B_n t$ instead, which is
smaller by a factor of order $\sqrt{n/\log n}$.

The bound in \cref{thm:probably-exact} is uniform over the base cost. The
Gaussian radius identity also quantifies how the base cost changes the radius
tail relative to the zero base cost. Define the \emph{assignment-cost range}
\begin{equation}\label{eq:assignment-cost-range}
  A_n(\Cbar) \coloneqq \max_{\pi \in \Sn} \ip{\Cbar, P_\pi} - \min_{\pi \in \Sn} \ip{\Cbar, P_\pi}.
\end{equation}

\begin{corollary}[Comparison with zero base cost]\label{cor:base-comparison}
  In the setting of \cref{thm:probably-exact}, for every $\Cbar \in \RR^{n
  \times n}$ and every $r \ge 0$,
  \begin{equation}\label{eq:base-universality}
    1 \le \frac{
      \PP[\rinf(\Cbar + \sigma G) \ge r]
    }{
      \PP[\rinf(\sigma G) \ge r]
    }
    \le e^{2 r A_n(\Cbar)/\sigma^2}.
  \end{equation}
\end{corollary}

\begin{proof}[Proof of \cref{cor:base-comparison}]
  Using the same Gaussian matrix $G$, let $\pi^\star \coloneqq
  \pi^\star(\Cbar + \sigma G)$ and $\pi_0 \coloneqq \pi^\star(\sigma G)$ be the
  optimal assignments with base cost $\Cbar$ and with zero base cost.
  \Cref{lem:margin-identity} with $u = r$, applied to both costs, gives
  \begin{align*}
    \PP[\rinf(\Cbar + \sigma G) \ge r]
    &= e^{-2(n-1) r^2/\sigma^2} \,
      \EE \exp\Bigl(\frac{2r}{\sigma} \ip{G, P_{\pi^\star} - J}\Bigr), \\
    \PP[\rinf(\sigma G) \ge r]
    &= e^{-2(n-1) r^2/\sigma^2} \,
      \EE \exp\Bigl(\frac{2r}{\sigma} \ip{G, P_{\pi_0} - J}\Bigr).
  \end{align*}
  The first exponential is the second times $\exp(\frac{2r}{\sigma} \ip{G,
  P_{\pi^\star} - P_{\pi_0}})$.  For every $G$, optimality of $\pi_0$ for
  $\sigma G$ and of $\pi^\star$ for $\Cbar + \sigma G$ implies
  \[
    0 \le \ip{G, P_{\pi^\star} - P_{\pi_0}}
    \le \frac{1}{\sigma} \ip{\Cbar, P_{\pi_0} - P_{\pi^\star}}
    \le \frac{1}{\sigma} A_n(\Cbar),
  \]
  so the ratio of the first exponential to the second lies between $1$ and
  $e^{2 r A_n(\Cbar)/\sigma^2}$. Both exponentials are positive, so their
  expectations satisfy the same bounds.  Canceling the common prefactor yields
  \cref{eq:base-universality}.
\end{proof}

The ratio in \cref{eq:base-universality} tends to one whenever $r
A_n(\Cbar)/\sigma^2 \to 0$. If $\Cbar = a \oplus b$ for some $a, b \in \RR^n$,
then every assignment has the same base cost, so $A_n(\Cbar) = 0$ and the two
tails in \cref{eq:base-universality} coincide.

\subsection{Random base costs and anisotropic noise}

\Cref{thm:probably-exact} is stated for a fixed base cost and isotropic noise.
Because the bound in \cref{eq:gap-law} does not depend on $\Cbar$, it extends to
random base costs and to anisotropic Gaussian noise.

\begin{corollary}[Random base costs and anisotropic noise]\label{cor:random-base}
  In the setting of \cref{thm:probably-exact}, let $\Cbar$ be random and
  independent of $G$.  Then the conclusions of the theorem hold.  They also
  hold for every Gaussian cost $C$ with nonsingular $\Gamma \coloneqq \Cov(C)$,
  where $\sigma \coloneqq \sqrt{\lambda_{\min}(\Gamma)} > 0$.
\end{corollary}

\begin{proof}[Proof of \cref{cor:random-base}]
  The assertions on the event $\{\rinf(C) \ge \sigma t\}$ are deterministic, so
  only the probability bound needs proof.  By \cref{lem:margin-forms}, the
  event is a finite union of translated optimality cones, so it is closed. Let
  $\nu$ be the law of $\Cbar$.  Because $\Cbar$ and $G$ are independent,
  Fubini's theorem lets us integrate over $G$ first:
  \[
    \PP[\rinf(C) \ge \sigma t]
    \;=\; \int \PP[\rinf(c + \sigma G) \ge \sigma t] \, \nu(dc).
  \]
  For each $c$, \cref{eq:gap-law} applied with base cost $c$ bounds the
  integrand below by a quantity that does not depend on $c$.  Because $\nu$ is
  a probability measure, the integral satisfies the same bound.

  For a Gaussian cost $C$, the event depends on $C$ alone, so its probability
  depends only on the law of $C$.  Since $\sigma^2 = \lambda_{\min}(\Gamma)$,
  the matrix $\Gamma - \sigma^2 I$ is positive semidefinite.  So there is a
  Gaussian $\Cbar$, possibly degenerate, independent of $G$, with the mean of
  $C$ and covariance $\Gamma - \sigma^2 I$.  Then $\Cbar + \sigma G$ is
  Gaussian with the mean and covariance of $C$, so it has the law of $C$, and
  the first part applies.
\end{proof}

Independent heteroscedastic noise is the diagonal case: if the entries of $C$
are independent Gaussian variables with standard deviations $s_{ij} > 0$, then
$\Gamma$ is diagonal and $\sigma = \min_{i, j} s_{ij}$. The probability bounds
in \cref{cor:fixed-confidence,cor:safe-scale} below follow from
\cref{eq:gap-law}, so both corollaries hold under either hypothesis of
\cref{cor:random-base}. Conditioning on $\Cbar$, after the same reduction in law
for Gaussian costs, also extends \cref{eq:beier-vocking-gap} under both
hypotheses of \cref{cor:random-base}, so the comparison with
\citet{BeierVocking2006} at the same $\sigma$ and failure probability is
unchanged.  Cost laws without an independent Gaussian component, such as
independent uniform costs, lie outside the corollary; see
\cref{sec:dictionary}.

\subsection{Fixed and vanishing failure probabilities}

\Cref{thm:probably-exact} holds for every threshold $t > 0$, and
\cref{eq:gap-law} bounds the probability of the event $\{\rinf(C) \ge \sigma
t\}$ from below by a function that decreases in $t$: increasing the threshold
raises the guaranteed radius but lowers the guaranteed confidence.  We apply
\cref{thm:probably-exact} with two choices of $t$. First,
\cref{cor:fixed-confidence} guarantees a radius of order $\sigma/(n \sqrt{\log
n})$ at any prescribed failure probability. Second, the failure bound vanishes
precisely when the threshold sequence satisfies $t_n = o(1/(n \sqrt{\log n}))$,
and in this case, \cref{cor:safe-scale} gives the explicit choice $t_n =
n^{-3/2}$. For both choices, write $t = x/(2 B_n)$. The resulting bound depends
on $n$ only through a term that tends to zero as $n \to \infty$. The section
concludes with a comparison to the bound~\cref{eq:beier-vocking-gap} of
\citet{BeierVocking2006}.

Write the right-hand side of the tail bound \cref{eq:gap-law} as $1 -
\delta_n(t)$, with failure bound
\begin{equation}\label{eq:failure-probability}
  \delta_n(t) \coloneqq 1
  - \exp(
    -2 B_n t
    - 2 (n-1) t^2
  ).
\end{equation}
For $x > 0$, set $t \coloneqq x/(2 B_n)$ to make the linear term $-2 B_n t$ in
the exponent equal to $-x$.  The radius level of \cref{thm:probably-exact} is
then $\sigma t = \sigma x/(2 B_n)$ and its gap level is $4 \sigma t = (2
\sigma/B_n) \cdot x$.  The identity $B_n^2 = 2(n-1) \log(n!)$ implies $2 (n-1)
t^2 = x^2/(4 \log(n!))$. Thus for every $x > 0$, \cref{eq:gap-law} reduces to
\begin{equation}\label{eq:scale-free-tail}
  \PP\Bigl[\Delta(C) \ge \frac{2 \sigma}{B_n} \cdot x\Bigr]
  \;\ge\; \PP\Bigl[\rinf(C) \ge \frac{\sigma}{2 B_n} \cdot x\Bigr]
  \;\ge\; \exp\Bigl(-x - \frac{x^2}{4 \log(n!)}\Bigr) .
\end{equation}
For each fixed $x > 0$, the lower bound in \cref{eq:scale-free-tail} increases
to $e^{-x}$ as $n$ grows.  Solving this bound for a prescribed failure
probability fixes the radius level.  For $\varepsilon \in (0, 1)$, set
$L_\varepsilon \coloneqq -\log(1 - \varepsilon)$, so that $e^{-L_\varepsilon} = 1 -
\varepsilon$, and let $x_\varepsilon$ be the positive root of $x + x^2/(4 \log(n!)) =
L_\varepsilon$.  Solving the quadratic gives the radius level
\begin{equation}\label{eq:fixed-confidence-margin}
  \bar{r}_\varepsilon(n, \sigma)
  \;\coloneqq\;
  \frac{\sigma x_\varepsilon}{2 B_n}
  \;=\;
  \frac{\sigma}{B_n} \cdot
  \frac{L_\varepsilon}{1 + \sqrt{1 + L_\varepsilon/\log(n!)}} \, .
\end{equation}

\begin{corollary}[Stability radius at prescribed confidence]\label{cor:fixed-confidence}
  In the setting of \cref{thm:probably-exact}, fix $\varepsilon \in (0, 1)$. Then,
  uniformly in $\Cbar$, with $C = \Cbar + \sigma G$,
  \[
    \PP[\rinf(C) \ge \bar{r}_\varepsilon(n, \sigma)]
    \;\ge\; 1 - \varepsilon,
  \]
  and on this event $\Delta(C) \ge 4 \bar{r}_\varepsilon(n, \sigma)$. Moreover,
  $\bar{r}_\varepsilon(n, \sigma) \sim \sigma L_\varepsilon/(2 B_n)$ as $n \to
  \infty$, which is of order $\sigma/(n \sqrt{\log n})$.
\end{corollary}

\begin{proof}[Proof of \cref{cor:fixed-confidence}]
  At $x = x_\varepsilon$, the right-hand side of \cref{eq:scale-free-tail} is
  $e^{-L_\varepsilon} = 1 - \varepsilon$, and the radius level there is
  $\bar{r}_\varepsilon(n, \sigma)$.  On this event, $\Delta(C) \ge 4 \rinf(C) \ge 4
  \bar{r}_\varepsilon(n, \sigma)$ by \cref{eq:cone-distance-bounds}. As $n \to
  \infty$, $\log(n!) \to \infty$, so the square root in
  \cref{eq:fixed-confidence-margin} tends to $1$ and $\bar{r}_\varepsilon(n,
  \sigma) \sim \sigma L_\varepsilon/(2 B_n)$. The order follows from $B_n \sim n
  \sqrt{2 \log n}$; see the remark after \cref{eq:Yn-Bn}.
\end{proof}

\Cref{sec:centered-base} shows that the radius level in
\cref{cor:fixed-confidence} is asymptotically sharp: additive base costs, which
are constant across all assignments, are the worst case.

By \cref{eq:failure-probability}, the negative logarithm $-\log(1 -
\delta_n(t_n))$ equals $2 t_n B_n + (t_n B_n)^2/\log(n!)$. Both terms are
nonnegative, so the failure bound $\delta_n(t_n)$ tends to zero if and only if
$t_n B_n \to 0$, that is, if and only if $t_n = o(1/(n \sqrt{\log n}))$. For
$t_n = n^{-a}$, this holds exactly when $a > 1$. The following corollary sets
$a = 3/2$ for definiteness.

\begin{corollary}[Stability radius with vanishing failure probability]\label{cor:safe-scale}
  In the setting of \cref{thm:probably-exact}, set $\bar{r}(n, \sigma)
  \coloneqq \sigma n^{-3/2}$ and $\delta(n) \coloneqq 2 \sqrt{2 \log n/n} +
  2/n^2$.  Then $\delta_n(n^{-3/2}) \le \delta(n)$, and
  \[
    \PP[\rinf(C) \ge \bar{r}(n, \sigma)]
    \;\ge\; 1 - \delta(n),
  \]
  uniformly in $\Cbar$, with $\delta(n) \to 0$.  On this event, $\Delta(C) \ge
  4 \bar{r}(n, \sigma)$.
\end{corollary}

\begin{proof}[Proof of \cref{cor:safe-scale}]
  Take $t = n^{-3/2}$, so that $\sigma t = \bar{r}(n, \sigma)$.  With $\log(n!)
  \le n \log n$, $n - 1 \le n$, and $1 - e^{-x} \le x$,
  \[
    \delta_n(n^{-3/2})
    \;\le\; 2 n^{-3/2} \sqrt{2 n \log(n!)}
    + 2 n^{-3}(n-1)
    \;\le\; 2 \sqrt{\frac{2 \log n}{n}} + \frac{2}{n^2}
    \;=\; \delta(n).
  \]
  Applying \cref{eq:gap-law} at $t = n^{-3/2}$ establishes the probability
  inequality, because $1 - \delta_n(n^{-3/2}) \ge 1 - \delta(n)$.  On this
  event, $\Delta(C) \ge 4 \rinf(C) \ge 4 \bar{r}(n, \sigma)$ by
  \cref{eq:cone-distance-bounds}.
\end{proof}

\begin{figure}[t]
  \centering
  \includegraphics{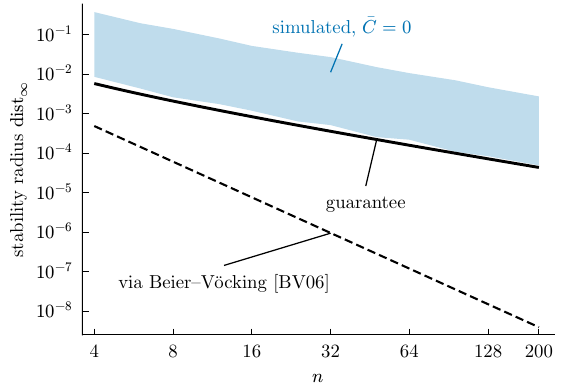}
  \caption{Stability radius against $n$ ($\sigma = 1$, $\varepsilon = 0.05$).
    Solid: the radius level $\bar{r}_\varepsilon(n, \sigma)$ guaranteed with
    probability $1 - \varepsilon$, uniformly in $\Cbar$, by
    \cref{cor:fixed-confidence}.  Shaded: 5th to 95th percentiles of
    $\rinf(C)$ at $\Cbar = 0$, the worst case (\cref{cor:extremal-base}), over
    $1000$ instances per size ($2000$ for $n \ge 128$).  Dashed: the radius
    guaranteed at the same $\varepsilon$ by \cref{eq:beier-vocking-gap} through
    $\rinf(C) \ge \Delta(C)/(2n)$.}
  \label{fig:tolerance}
\end{figure}

At a common failure probability $\varepsilon$, the Beier--Vöcking bound
\cref{eq:beier-vocking-gap} guarantees the gap $\varepsilon \sigma \sqrt{2
\pi}/(2 n^2)$ and, through $\rinf(C) \ge \Delta(C)/(2n)$, the radius $\varepsilon
\sigma \sqrt{2 \pi}/(4 n^3)$. On the other hand, our \cref{cor:fixed-confidence}
guarantees the radius $\bar{r}_\varepsilon(n, \sigma)$ and the gap $4
\bar{r}_\varepsilon(n, \sigma)$.  The ratios of our guarantees to theirs are
therefore asymptotic to
\begin{equation}\label{eq:asymptotic-ratios}
  \frac{L_\varepsilon}{\varepsilon} \cdot 2 \sqrt{2/\pi} \, \frac{n^2}{B_n}
  \quad\text{for the gap,}\quad\text{and}\quad
  \frac{L_\varepsilon}{\varepsilon} \cdot \sqrt{2/\pi} \, \frac{n^3}{B_n}
  \quad\text{for the radius},
\end{equation}
respectively, of orders $n/\sqrt{\log n}$ and $n^2/\sqrt{\log n}$.  The radius
factor is $n/2$ times the gap factor because the two bounds pass between gap
and radius at different rates: ours through $\Delta(C) \ge 4 \rinf(C)$, theirs
through $\rinf(C) \ge \Delta(C)/(2n)$.  Since $L_\varepsilon \ge \varepsilon$ and
$L_\varepsilon/\varepsilon \to 1$ as $\varepsilon \to 0$, each factor in
\cref{eq:asymptotic-ratios} is at least its limit as $\varepsilon \to 0$.
\Cref{fig:tolerance} compares, at $\varepsilon = 0.05$, the stability radius
guaranteed by \cref{cor:fixed-confidence} with the radius obtained from
\cref{eq:beier-vocking-gap} and with simulated radii at $\Cbar = 0$.

\section{Extremal base costs and the exponential limit}\label{sec:centered-base}

\begin{figure}[t]
  \centering
  \includegraphics{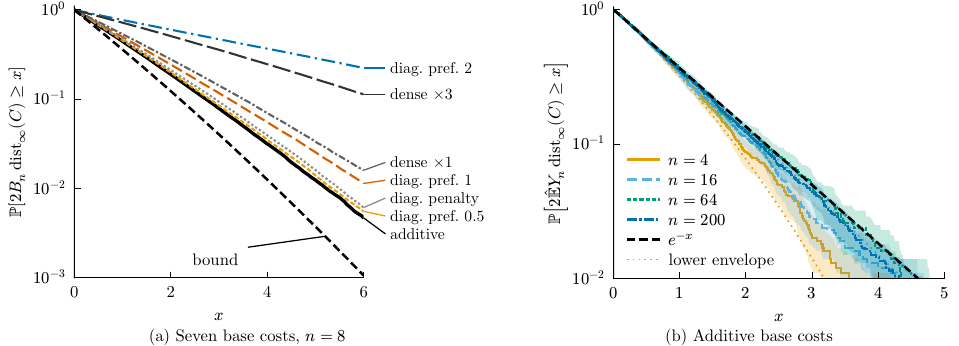}
  \caption[Empirical survival functions of the stability radius]{Empirical
    survival functions of the stability radius ($\sigma = 1$).
    (a) Radii scaled by $2 B_n$, with $n = 8$, over $10^5$ Gaussian draws
    shared by seven base costs: $\Cbar = 0$ (additive); $\lambda(\ones
    \ones^\top - I)$ with $\lambda = 0.5, 1, 2$ (diagonal preference); $10 I$
    (diagonal penalty); and $s Z$ with $s = 1, 3$, where $Z$ is one fixed
    standard Gaussian matrix (dense). Dashed: the lower bound in
    \cref{eq:scale-free-tail}.
    (b) Additive base costs, $\Cbar = 0$. Scaling by $2 \EE Y_n$ makes the
    theoretical upper envelope $e^{-x}$ for every $n$; see
    \cref{eq:centered-margin-sandwich}. The plots substitute the sample mean
    $\hat{\EE} Y_n$ in the normalization and lower envelopes, using the same
    instances ($1000$ per size, $2000$ at $n = 200$). Dashed: $e^{-x}$, also
    the limiting survival function in \cref{cor:extremal-base}(ii); dotted:
    estimated lower envelopes from \cref{eq:centered-margin-sandwich}; bands:
    pointwise $95\%$ bootstrap intervals. Radii are checked against enumeration
    for $3 \le n \le 7$.}
  \label{fig:margin}
\end{figure}

At a fixed confidence level, \cref{cor:fixed-confidence} guarantees a stability
radius of order $\sigma/(n \sqrt{\log n})$, uniformly over the base cost
$\Cbar$. Moreover, this uniform guarantee is asymptotically sharp: additive
base costs (including $\Cbar = 0$) are the worst case, and at these base costs
the scaled radius has an exponential limit (\cref{fig:margin}).

Suppose the base cost $\Cbar$ is \emph{additive}, that is, $\Cbar = a \oplus b$
for some $a, b \in \RR^n$. Then the inner product $\ip{\Cbar, P_\pi} = \ip{a,
\ones} + \ip{b, \ones}$ is independent of $\pi$. Thus, $\pi^\star(\Cbar + \sigma
G) = \pi^\star(\sigma G)$ and, by \cref{eq:dist-p-def}, $\rinf(\Cbar + \sigma G)
= \rinf(\sigma G)$. The following corollary identifies additive base costs as
exactly those at which the stability radius is stochastically minimal. It also
gives the asymptotic rate of the radius tail and the limit law of the scaled
radius $(2 B_n/\sigma) \, \rinf(\sigma G)$, with $B_n$ from \cref{eq:Yn-Bn}.

\begin{corollary}[Extremal base costs and exponential limit]\label{cor:extremal-base}
  In the setting of \cref{thm:probably-exact}:
  \begin{enumerate}[(i)]
    \item For every $r > 0$,
      \[
        \min_{\Cbar \in \RR^{n \times n}}
        \PP[\rinf(\Cbar + \sigma G) \ge r]
        \;=\; \PP[\rinf(\sigma G) \ge r].
      \]
      The minimum is attained exactly at the additive base costs.
    \item Suppose the base costs are additive.  As $n \to \infty$, the
      logarithmic tail satisfies
      \begin{equation}\label{eq:centered-margin-rate}
        -\log \PP[\rinf(C) \ge \sigma t_n]
        \;\sim\; 2 t_n B_n
      \end{equation}
      for every positive sequence $t_n = o(\sqrt{\log n})$. In particular,
      \begin{equation}\label{eq:centered-margin-limit}
        \sup_{x \ge 0} \bigl\lvert \PP[(2 B_n/\sigma) \rinf(C) \ge x] - e^{-x}\bigr\rvert
        \;\to\; 0,
      \end{equation}
      so $(2 B_n/\sigma) \rinf(C)$ converges in distribution to a standard
      exponential variable.  Moreover, $\EE \rinf(C) \sim \sigma/(2 B_n)$.
  \end{enumerate}
\end{corollary}

\begin{proof}[Proof of \cref{cor:extremal-base}]
  (i) Fix $r > 0$, set $t \coloneqq r/\sigma$ and $q_n(t) \coloneqq 2 (n-1)
  t^2$, and write $H_\pi \coloneqq \ip{G, P_\pi - J}$ as in the proof
  of~\cref{thm:probably-exact}. By \cref{cor:base-comparison}, the minimum
  equals $\PP[\rinf(\sigma G) \ge r]$, and every additive base cost attains it.
  For additive base costs, the permutation $\pi^\star(C)$ minimizes $\ip{G,
  P_\pi}$, so $H_{\pi^\star(C)} = -Y_n$. The Gaussian radius identity
  \cref{eq:smoothed-mixture} therefore reduces to
  \begin{equation}\label{eq:smoothed-mixture-centered}
    \PP[\rinf(\sigma G) \ge \sigma t]
    \;=\; e^{-q_n(t)} \,
    \EE e^{-2 t Y_n}.
  \end{equation}

  Conversely, suppose that $\Cbar$ is not additive. Set $a_i \coloneqq
  \Cbar_{i1} - \Cbar_{11}$ and $b_k \coloneqq \Cbar_{1k}$ for all $i$ and $k$.
  Set $D \coloneqq \Cbar - a \oplus b$. Its first row and first column are
  zero, but $D \ne 0$ because $\Cbar$ is not additive. Thus $D_{ik} \ne 0$ for
  some $i \ne 1$ and $k \ne 1$. Choose $\pi$ with $\pi(i) = k$ and $\pi(1) = 1$,
  and let $\pi'$ agree with $\pi$ except that $\pi'(i) = 1$ and $\pi'(1) = k$.
  Then
  \[
    \ip{\Cbar, P_\pi - P_{\pi'}}
    = \Cbar_{ik} + \Cbar_{11} - \Cbar_{i1} - \Cbar_{1k}
    = D_{ik} \ne 0 .
  \]
  Relabeling $\pi$ and $\pi'$ if needed, we obtain $\rho \ne \pi$ with $d
  \coloneqq \ip{\Cbar, P_\pi - P_\rho} > 0$. The open set
  \[
    \mathcal{U}_\pi \coloneqq \{Z : \ip{Z, P_\tau - P_\pi} > 0 \text{ for all } \tau \ne \pi,
      \; \ip{Z, P_\rho - P_\pi} < d/\sigma\}
  \]
  contains $-\lambda P_\pi$ for $0 < \lambda < d/(2 \sigma)$, and hence
  $\PP[G \in \mathcal{U}_\pi] > 0$. On the event $G \in \mathcal{U}_\pi$, the
  permutation $\pi$ uniquely minimizes $\ip{G, P_\tau}$ over $\tau$, so $-Y_n =
  H_\pi$.  But $\ip{\Cbar + \sigma G, P_\rho - P_\pi} < 0$, so $\pi^\star(C) \ne
  \pi$ and $H_{\pi^\star(C)} > -Y_n$. The difference $e^{2 t H_{\pi^\star(C)}}
  - e^{-2 t Y_n}$ is nonnegative everywhere and positive on $\mathcal{U}_\pi$,
  an event of positive probability. Hence, by \cref{lem:margin-identity} and
  \cref{eq:smoothed-mixture-centered}, the lower inequality in
  \cref{eq:base-universality} is strict, and this non-additive base cost does
  not attain the minimum.

  (ii) Write $S_n(t) \coloneqq \PP[\rinf(C) \ge \sigma t]$ for the tail
  probability in~\cref{eq:smoothed-mixture-centered}.  We first bound $S_n(t)$
  above and below by exponentials in $\EE Y_n$; the limit $\EE Y_n/B_n \to 1$
  \citep[Thm.~2.3]{MordantSegers2021} then establishes the tail rate and limit
  law.

  Jensen's inequality, as in~\cref{eq:jensen}, applied to the right-hand side
  of~\cref{eq:smoothed-mixture-centered} yields the lower bound $e^{-q_n(t) - 2
  t \EE Y_n} \le S_n(t)$. For the upper bound, factor the expectation on that
  right-hand side as $\EE e^{-2 t Y_n} = e^{-2 t \EE Y_n} \, \EE e^{-2 t (Y_n -
  \EE Y_n)}$. As a function of $G$, $Y_n$ is the maximum of the linear
  functions $\ip{G, J - P_\pi}$ of~\cref{eq:Yn-Bn}.  Their gradients have
  Frobenius norm $\sqrt{n-1}$ by~\cref{eq:hamming-identities}, so $Y_n$ is
  $\sqrt{n-1}$-Lipschitz in the Frobenius norm. The Gaussian exponential-moment
  bound \citep[Thm.~5.5]{BoucheronLugosiMassart2013} at parameter $-2 t$ then
  implies
  \[
    \EE e^{-2 t (Y_n - \EE Y_n)}
    \;\le\; e^{2 (n-1) t^2}
    = e^{q_n(t)},
  \]
  and this factor cancels the prefactor $e^{-q_n(t)}$
  in~\cref{eq:smoothed-mixture-centered}. Combining the upper and lower bounds
  yields
  \begin{equation}\label{eq:centered-margin-sandwich}
    e^{-q_n(t) - 2 t \EE Y_n}
    \;\le\;
    S_n(t)
    \;\le\;
    e^{-2 t \EE Y_n}.
  \end{equation}
  Under the normalization $x = 2 t \EE Y_n$, the upper bound in
  \cref{eq:centered-margin-sandwich} becomes $e^{-x}$ and the lower bound
  becomes $\exp(-x - (n-1) x^2/(2 (\EE Y_n)^2))$.

  For \cref{eq:centered-margin-rate}, take negative logarithms in
  \cref{eq:centered-margin-sandwich} at $t = t_n$, and divide by $2 t_n B_n$:
  \[
    \frac{\EE Y_n}{B_n}
    \;\le\; -(2 t_n B_n)^{-1} \log S_n(t_n)
    \;\le\; \frac{\EE Y_n}{B_n} + \frac{q_n(t_n)}{2 t_n B_n} .
  \]
  The hypothesis $t_n = o(\sqrt{\log n})$ implies $q_n(t_n)/(2 t_n B_n) = t_n
  (n-1)/B_n = O(t_n/\sqrt{\log n}) \to 0$, and $\EE Y_n/B_n \to 1$.  Both sides
  therefore tend to $1$.

  For the limit law, fix $x > 0$ and take $t_n = x/(2 B_n) = o(\sqrt{\log n})$
  in \cref{eq:centered-margin-rate}: its right-hand side equals $x$, so $-\log
  S_n(x/(2 B_n)) \to x$ and hence $S_n(x/(2 B_n)) \to e^{-x}$; at $x = 0$ both
  survival functions equal one. Since $t \mapsto 2 t B_n$ maps $[0, \infty)$
  onto itself, the supremum in \cref{eq:centered-margin-limit} is the uniform
  distance between the survival function of $(2 B_n/\sigma) \rinf(C)$ and
  $e^{-x}$. By Pólya's theorem \citep[Lem.~2.11]{vanderVaart1998}, pointwise
  convergence of survival functions to a continuous survival function is
  uniform; this proves \cref{eq:centered-margin-limit} and the convergence in
  distribution.

  For convergence of the mean, at $t = x/(2 B_n)$ the upper bound in
  \cref{eq:centered-margin-sandwich} is $e^{-x \EE Y_n/B_n} \le e^{-\beta x}$,
  where $\beta \coloneqq \inf_{n \ge 2} \EE Y_n/B_n > 0$ because each ratio is
  positive and the ratios tend to one.  Integrating the survival functions and
  applying dominated convergence yields $\EE \rinf(C) \sim \sigma/(2 B_n)$.
\end{proof}

\Cref{fig:margin}(b) compares the empirical tail $\PP[\rinf(\sigma G) \ge \sigma
t]$ with the two bounds in \cref{eq:centered-margin-sandwich}.

\Cref{cor:extremal-base} shows that the radius bound in
\cref{eq:scale-free-tail} is asymptotically sharp after minimization over base
costs. For each fixed $x > 0$, part (i) says that the minimum over $\Cbar$ of
$\PP[\rinf(\Cbar + \sigma G) \ge \sigma x/(2 B_n)]$ is attained at additive base
costs, and part (ii) says that this minimum tends to $e^{-x}$. Thus, for every
fixed $\varepsilon \in (0, 1)$ and every fixed $\alpha > 0$, the corresponding
limit at the radius level of \cref{cor:fixed-confidence} is
\begin{equation}\label{eq:fixed-confidence-sharpness}
  \lim_{n \to \infty} \min_{\Cbar \in \RR^{n \times n}} \PP[\rinf(\Cbar + \sigma G) \ge \alpha \bar{r}_\varepsilon(n, \sigma)]
  \;=\; (1 - \varepsilon)^\alpha.
\end{equation}
Indeed, the scaled radius level $(2 B_n/\sigma) \bar{r}_\varepsilon(n, \sigma)$
tends to $L_\varepsilon$ by \cref{cor:fixed-confidence}, so the uniform convergence
in \cref{eq:centered-margin-limit} gives the limiting minimum in
\cref{eq:fixed-confidence-sharpness}. For every fixed $\alpha > 1$, the limiting
probability $(1 - \varepsilon)^\alpha$ is below $1 - \varepsilon$. Thus the radius
level $\bar{r}_\varepsilon(n, \sigma)$ cannot be increased by a fixed factor
greater than one while retaining confidence $1 - \varepsilon$ uniformly over base
costs for all sufficiently large $n$.

Only additive base costs attain the minimum in \cref{cor:extremal-base}(i). The
exponential limit and mean asymptotic, however, hold more broadly. Recall that
the assignment-cost range $A_n(\Cbar)$, defined in
\cref{eq:assignment-cost-range}, is the difference between the largest and
smallest assignment costs under $\Cbar$. It vanishes precisely when the base
costs are additive. The comparison in \cref{cor:base-comparison} yields the
following extension to sequences of base costs whose range is $o(\sigma B_n)$.

\begin{corollary}[Base costs near additive]\label{cor:near-additive}
  In the setting of \cref{thm:probably-exact}, if a sequence of base costs
  $\Cbar_n$ satisfies $A_n(\Cbar_n) = o(\sigma B_n)$, then the scaled radius $(2
  B_n/\sigma) \rinf(\Cbar_n + \sigma G)$ converges in distribution to a
  standard exponential variable, and $\EE \rinf(\Cbar_n + \sigma G) \sim
  \sigma/(2 B_n)$. The hypothesis holds if $\inf_{a, b \in \RR^n}
  \infnorm{\Cbar_n - a \oplus b} = o(\sigma \sqrt{\log n})$, and hence, for
  fixed $\sigma$, for every sequence of base costs with uniformly bounded
  entries.
\end{corollary}

\begin{proof}[Proof of \cref{cor:near-additive}]
  By \cref{cor:base-comparison}, for every fixed $x \ge 0$, at $r = \sigma x/(2
  B_n)$ the factor in \cref{eq:base-universality} is $e^{x A_n(\Cbar_n)/(\sigma
  B_n)}$, which tends to one, so \cref{eq:centered-margin-limit} implies the
  distributional limit. Combining \cref{eq:base-universality} with
  \cref{eq:centered-margin-sandwich} bounds the normalized survival function
  above by
  \[
    \exp\Bigl[-x \Bigl(\frac{\EE Y_n}{B_n} - \frac{A_n(\Cbar_n)}{\sigma B_n}\Bigr)\Bigr]
    \le e^{-x/2},
  \]
  for all $x \ge 0$ and all sufficiently large $n$, because $\EE Y_n/B_n \to 1$
  and $A_n(\Cbar_n)/(\sigma B_n) \to 0$. Dominated convergence therefore
  establishes the mean asymptotic, as in the proof of
  \cref{cor:extremal-base}(ii). Finally, subtracting an additive matrix leaves
  $A_n$ unchanged and each assignment uses $n$ entries, so $A_n(\Cbar_n) \le 2 n
  \inf_{a, b \in \RR^n} \infnorm{\Cbar_n - a \oplus b}$. Together with $B_n \sim
  n \sqrt{2 \log n}$, this bound shows that the residual-norm condition
  suffices and that uniformly bounded entries are covered.
\end{proof}

\section{Discussion}\label{sec:dictionary}

\begin{figure}[t]
  \centering
  \includegraphics{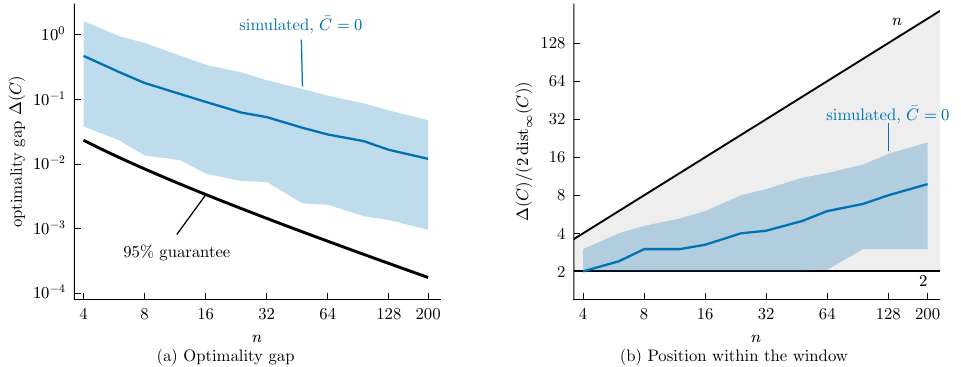}
  \caption{The optimality gap at $\Cbar = 0$ and $\sigma = 1$, over the
    instances of \cref{fig:tolerance}. Blue: median (line) and 5th to 95th
    percentiles (band). (a) The gap $\Delta(C)$. Black: the level $4
    \bar{r}_\varepsilon(n, \sigma)$, with $\varepsilon = 0.05$, guaranteed with
    probability at least $95\%$, uniformly over base costs, by
    \cref{cor:fixed-confidence}. (b) The ratio $\Delta(C)/(2 \rinf(C))$. Gray:
    the interval $[2, n]$ given by the gap--radius bounds
    \cref{eq:cone-distance-bounds}.}
  \label{fig:gap}
\end{figure}

The stability radius bound in \cref{thm:probably-exact} is uniform over base
costs and asymptotically sharp. The optimality gap, however, is less well
understood. Even at additive base costs, where \cref{cor:extremal-base} gives
the distributional limit of the scaled radius $R_n \coloneqq (2 B_n/\sigma)
\rinf(C)$, we still lack a precise asymptotic characterization. In particular,
rewriting the gap--radius bounds \cref{eq:cone-distance-bounds} in terms of
$R_n$ gives
\[
  2 \frac{\sigma R_n}{B_n} \le \Delta(C) \le n \frac{\sigma R_n}{B_n},
\]
which leaves a factor of $\half n$ between the upper and lower bounds.
\Cref{fig:gap}(a) shows the simulated gap at $\Cbar = 0$ beside the gap level $4
\bar{r}_\varepsilon(n, \sigma)$ guaranteed by \cref{cor:fixed-confidence}.
\Cref{fig:gap}(b) plots the ratio $\Delta(C)/(2 \rinf(C))$ within its bounds
$[2, n]$. Over the sampled sizes, the empirical median increases with $n$ but
stays far below $n$. The dependence of the gap on the base costs raises a
related question. Additive base costs minimize the radius tail by
\cref{cor:extremal-base}(i), but do they also minimize $\PP[\Delta(\Cbar +
\sigma G) \ge s]$ over $\Cbar$ for each $s > 0$?

The extension to non-Gaussian noise in \cref{cor:random-base} still requires an
independent Gaussian component $\sigma G$. Independent uniform noise provides a
concrete model without such a component. Let $C = \Cbar + \sigma U$, where $U$
has independent uniform $[0, 1]$ entries. The generalized isolating lemma of
\citet[Lem.~5]{BeierVocking2006}, combined with \cref{eq:cone-distance-bounds},
guarantees a radius of $\varepsilon \sigma/(4 n^3)$ with probability at least $1 -
\varepsilon$, uniformly in $\Cbar$. The Gaussian density calculation in
\cref{eq:shifted-cone-measure} does not apply to this model. Can the
uniform-noise guarantee be improved to a larger order, perhaps to the Gaussian
scale $\sigma/(n \sqrt{\log n})$?

The translated-cone argument also extends beyond the assignment problem to
other 0/1 polytopes. The proof of \cref{thm:probably-exact} relies on two
properties of permutation matrices: equidistance from the barycenter $J =
\frac{1}{n} \ones \ones^\top$ and the identity $\ip{J - P_\pi, P_\rho - P_\pi} =
\half \norm{P_\rho - P_\pi}_1$, which follows from
\cref{eq:hamming-identities}. Every polytope with vertices in $\{0, 1\}^d$
satisfies both properties with $J$ replaced by the projection of $\half \ones$
onto its affine hull. The radius event is then a disjoint union of translated
normal cones. Summing their Gaussian measures as in
\cref{eq:shifted-cone-measure} and bounding the integrand below by its minimum
over the vertices yields a tail bound uniform over base costs, with the
Gaussian width of the vertex set replacing $\EE Y_n$. At $\Cbar = 0$,
\citet[\S2.1.1]{FKPS2025} use instead a radial shift for inscribed polytopes.
Whether this extension of \cref{thm:probably-exact} is asymptotically sharp for
spanning trees, matroid bases, or perfect matchings of general graphs remains
open.

\section*{Declarations}

This paper continues my work with the late Paul Tseng \cite{FriedlanderTseng2007}
on exact regularization of convex programs and, more recently, on its average
case with Kubal, Plan, and Scott \cite{FKPS2025}. Here I treat the smoothed case
for the assignment problem. I used Claude (Anthropic), through Claude Code, and
GPT models (OpenAI), through Codex, to propose and test proof strategies, to
look for errors in proofs, to write and run the code for the computations in the
figures, and to copyedit and proofread the manuscript. I verified the proofs,
read the cited works, and tested the code. I assume full responsibility for all
content.

\bibliographystyle{plainnat-alpha}
\bibliography{refs}

\end{document}